%% file: control_slemma.tex
\documentclass[a4paper,12pt]{article}
\input defs.tex

\newcommand{\ackERC}{This project has received funding from the European Research Council (ERC) under the European Union's Horizon 2020 research and innovation programme grant agreement OCAL, No.\ 787845.}

\title{\headingfont Improving Tractability of Real-Time Control Schemes via Simplified \SL-Lemma}

\author{Goran Banjac, Jianzhe Zhen, Dick den Hertog, and John Lygeros}

\begin{document}

\maketitle

\begin{abstract}
  Various control schemes rely on a solution of a convex optimization problem involving a particular robust quadratic constraint, which can be reformulated as a linear matrix inequality using the well-known \SL-lemma.
  However, the computational effort required to solve the resulting semidefinite program may be prohibitively large for real-time applications requiring a repeated solution of such a problem.
  We use some recent advances in robust optimization that allow us to reformulate such a robust constraint as a set of linear and second-order cone constraints, which are computationally better suited to real-time applications.
  A numerical example demonstrates a huge speedup that can be obtained using the proposed reformulation.
\end{abstract}

\section{Introduction}\label{sec:intro}
Controlled dynamical systems are inevitably affected by unknown disturbances such as process noise, measurement noise, and imprecision in controller implementation.
One way to deal with such disturbances in optimal control problems is via \emph{stochastic optimization} in which optimization problems contain expectations and/or chance constraints \cite{Mesbah:2016,Heirung:2018}.
An important assumption in stochastic optimization is that the true joint probability distribution of the uncertain parameters is known or estimated, but in practice estimating such a distribution is  data-intensive \cite{Yanikoglu:2019}.
Moreover, the resulting optimization problems are computationally demanding as existing solution schemes often rely on sampling techniques \cite{Birge:2011}.

\emph{Robust optimization} is a powerful alternative approach to treat data uncertainty, where one seeks a control input that minimizes the worst-case cost across all possible realizations of disturbances included within a prescribed set.
Advantages of robust optimization include modeling flexibility, computational tractability, and inherent probabilistic guarantees associated to the obtained solutions \cite{Ben-Tal:2009}.

However, robust optimization is often criticized to be conservative as it
(i) focuses on minimizing the cost function with respect to the worst-case realization of disturbances, and
(ii) only takes into account the support of disturbances and ignores the associated statistical information, such as mean and variance, which are readily available or easy to estimate from the historical data.
To deal with the first issue, \emph{robust regret-optimal optimization} minimizes instead the worst-case regret, \ie the worst-case cost compared to the minimal cost given that the disturbance were known \emph{a priori} \cite{Delage:2019}.
\emph{Distributionally robust optimization} enables incorporation of statistical knowledge about the disturbance and minimizes the cost function with respect to its worst-case probability distribution \cite{Delage:2010,Wiesemann:2014,Chen:2020}.

Despite many appealing properties of the aforementioned techniques, it is essential to be able to solve the corresponding optimization problems efficiently.
When considering linear dynamical systems and quadratic cost functions, also known as \emph{linear quadratic control (LQC)}, many of these techniques rely on a solution of a convex optimization problem involving a particular robust quadratic constraint, which can be reformulated as a linear matrix inequality (LMI) using the celebrated \SL-lemma \cite{Yakubovich:1971}.
For instance, the authors in \cite{Bertsimas:2007} consider a robust optimization approach to constrained LQC under disturbances that are confined in an ellipsoidal set, and reformulate the resulting optimization problem as a semidefinite program (SDP).
Despite some recent advances in numerical methods for SDPs \cite{Majumdar:2020,Zheng:2020}, the computational effort required to solve these problems may be prohibitively large for real-time applications requiring a repeated solution of such a problem.

As an alternative, we propose using a simplified \SL-lemma, introduced in \cite{Ben-Tal:2014}, that reformulates such a robust quadratic constraint as a set of linear and second-order cone (SOC) constraints, which are computationally better suited to real-time applications than an LMI \cite{Domahidi:2013}.
Using this simplification, we reduce the SDP reformulation of the robust LQC problem proposed in \cite{Bertsimas:2007} to a second-order cone program (SOCP) and demonstrate on a numerical example that such an SOCP is up to three orders of magnitude faster than the SDP.
We also propose SOCP reformulations of other solution techniques, including robust regret-optimal control and distributionally robust control.

Moreover, we show that a simplified \SL-lemma is useful in other control schemes, including model predictive control (MPC).
Adding flexibility to an ellipsoidal terminal set in MPC, by allowing it to be scaled and translated online, is known to improve its performance.
We show that the resulting optimization problem can still be formulated as an SOCP.

\subsection*{Notation}
Let $\Nat$ denote the set of positive integers, $\Re$ the set of real numbers, $\Re_+$ the set of nonnegative real numbers, $\Re^n$ the $n$-dimensional Euclidean space, $\Re^{m\times n}$ the set of real $m$-by-$n$ matrices, and $\symm^n$ ($\symm_{++}^n$) the set of real $n$-by-$n$ symmetric (positive definite) matrices.
We denote the set of positive integers lower than or equal to $n\in\Nat$ by $[n]$, the vector of all ones by $\boldsymbol{1}$, the $i$-th component of $x\in\Re^n$ by $x_i$, the vertical concatenation of $x\in\Re^n$ and $y\in\Re^m$ by $(x,y)$, the operator mapping a vector to a diagonal matrix by $\diag\colon\Re^n\to\symm^n$, and $A\in\symm^n$ being positive (semi)definite by $A\succ 0$ ($A\succeq 0$).
For a vector $x\in\Re^n$, the inequality $x\ge 0$ is understood elementwise.

\section{Preliminaries}\label{sec:prelim}

\subsection{Second-Order Cone Programming}\label{subsec:socp}
An SOCP is a convex optimization problem given by \cite{Boyd:2004}
\[
  \begin{array}{ll}
    \underset{x}{\rm minimize} & f^T x \\
    \text{subject to} & \norm{A_i x + b_i}_2 \le c_i^T x + d_i, \quad i\in[m] \\
                      & Fx=g,
  \end{array}
\]
where $x\in\Re^n$ is the optimization variable, $A_i\in\Re^{n_i\times n}$, and~$F\in\Re^{p\times n}$.
We refer to the inequality constraints as SOC constraints and to $n_i$ as the dimension of the $i$-th SOC constraint.
We recall below a well-known SOC reformulation technique.
\begin{lemma}[\hspace{1sp}{\cite[\S 2.3]{Lobo:1998}}]\label{lem:hyperbolic}~
\begin{enumerate}[label=(\roman*)]
  \item $(y_i,z_i)\ge 0, \; x_i^2 \le y_i z_i \; \Longleftrightarrow \; \norm{(2x_i, y_i-z_i)}_2 \le y_i+z_i$.
  \item $a^Tx+b>0, \ \norm{Fx+g}_2^2 / (a^Tx+b) \le t$ \\[.2em] $\Longleftrightarrow \; \left\| \big( 2(Fx+g), t-a^Tx-b \big) \right\|_2 \le t+a^Tx+b$.
\end{enumerate}
\end{lemma}
We will make an extensive use of the first result in \Lem~\ref{lem:hyperbolic}, and thus define
\[
  \mcf{Q} \eqdef \left\lbrace (x_i,y_i,z_i)\in\Re^3 \mid \left\| (2x_i,y_i-z_i) \right\|_2 \le y_i+z_i \right\rbrace
\]
for ease of exhibition.
Observe that $\mcf{Q}$ is a cone, \ie $(x_i,y_i,z_i)\in\mcf{Q}$ implies $(rx_i,ry_i,rz_i)\in\mcf{Q}$ for all $r\ge 0$.

\subsection{A Simplified \SL-Lemma}\label{subsec:s-lemma}
Let $A\in\symm^n$, $D\in\symm^n$, $b\in\Re^n$, $e\in\Re^n$, $c\in\Re$, and $f\in\Re$.
We are interested in the following robust (not necessarily convex) quadratic constraint
\begin{equation}\label{eqn:quad-semi-inf}
  z^T D z + 2e^T z + f \ge 0, \quad \forall z \colon z^T A z + 2b^T z + c \ge 0,
\end{equation}
which is equivalent to the implication
\[
  z^T A z + 2b^T z + c \ge 0 \quad\Longrightarrow\quad z^T D z + 2e^T z + f \ge 0.
\]
Provided that there exists some $z_0\in\Re^n$ such that
\begin{equation}\label{eqn:strict}
  z_0^T A z_0 + 2b^T z_0 + c > 0,
\end{equation}
the classical \SL-lemma \cite{Yakubovich:1971} states that the implication above is valid if and only if there exists some $\lambda\in\Re_+$ such that
\begin{equation}\label{eqn:lmi}
  \left[ \begin{array}{c|c} D - \lambda A & e-\lambda b \\ \hline (e-\lambda b)^T & f - \lambda c \end{array} \right] \succeq 0.
\end{equation}
If $A$ and $D$ are \emph{simultaneously diagonalizable}, \ie there exists a nonsingular matrix $S\in\Re^{n\times n}$ such that
\begin{equation}\label{eqn:S}
  S^T A S = \diag(\alpha_1,\ldots,\alpha_n)
  \qquad\text{and}\qquad
  S^T D S = \diag(\delta_1,\ldots,\delta_n),
\end{equation}
then the LMI in \eqref{eqn:lmi} can be reduced to a set of linear and SOC constraints.

\begin{lemma}[\hspace{1sp}{\cite[\Thm~9]{Ben-Tal:2014}}]\label{lem:simplified-s}
Suppose there exist a vector $z_0\in\Re^n$ and a nonsingular matrix $S\in\Re^{n\times n}$ such that \eqref{eqn:strict} and \eqref{eqn:S} hold.
Then the robust constraint in \eqref{eqn:quad-semi-inf} holds if and only if there exists $(\lambda,t)\in\Re_+\times\Re^n$ such that
\begin{align*}
  & f - \lambda c \ge \boldsymbol{1}^T t \\
  & \big( \varepsilon_i - \lambda \beta_i, t_i, \delta_i-\lambda \alpha_i \big) \in \mcf{Q}, \quad \forall i\in[n],
\end{align*}
where $\varepsilon\eqdef S^T e$ and $\beta\eqdef S^T b$.
\end{lemma}

To shed some light on the constraints in \Lem~\ref{lem:simplified-s}, suppose for simplicity that $D-\lambda A \succ 0$, or equivalently
\[
  \delta_i - \lambda\alpha_i > 0, \quad i\in[n].
\]
Then \eqref{eqn:lmi} is equivalent to
\[
  \left[ \begin{array}{c|c} \diag(\delta_1,\ldots,\delta_n) - \lambda \diag(\alpha_1,\ldots,\alpha_n) & \varepsilon -\lambda \beta \\ \hline (\varepsilon -\lambda \beta)^T & f - \lambda c \end{array} \right] \succeq 0.
\]
Using the Schur complement \cite[\S A.5.5]{Boyd:2004}, the LMI reduces to
\[
  f - \lambda c \ge \sum_{i=1}^n \frac{(\varepsilon_i - \lambda \beta_i)^2}{\delta_i - \lambda\alpha_i},
\]
which, due to \Lem~\ref{lem:hyperbolic}, is equivalent to the constraints in \Lem~\ref{lem:simplified-s}.

\begin{remark}\label{rem:pos-def}
Positive definiteness of either $A$ or $D$ is a sufficient condition for the simultaneous diagonalizability of $A$ and $D$ \cite{Ben-Tal:2014}.
To show this, assume $A\succ 0$ and let $A=L L^T$ be its Cholesky decomposition.
Now define $C=L\inv D L^{-T}$ and let its orthogonal decomposition be given by $C=Q\diag(\delta_1,\ldots,\delta_n)Q^T$.
Then for $S=L^{-T} Q$, we have $S^T A S = I$ and $S^T D S=\diag(\delta_1,\ldots,\delta_n)$.
We refer the interested reader to \cite{Jiang:2016} for necessary and sufficient conditions for simultaneous diagonalizability of two matrices.
\end{remark}

\begin{remark}
A limitation of \Lem~\ref{lem:simplified-s} with respect to \eqref{eqn:lmi} is that $D$ cannot be an arbitrary affine map of optimization variables.
Also, in order to use an approximate \SL-lemma, which involves the intersection of multiple sets of the form $\lbrace z\in\Re^n \mid z^T A_k z + b_k^T z + c_k \ge 0 \rbrace$ for $k=1,\ldots,m$ \cite[\S 2.6.3]{Boyd:1994}, one would require that $D$ and all $A_k$ are jointly simultaneously diagonalizable.
\end{remark}

\section{Finite-Horizon LQC}\label{sec:lqr}
We consider an uncertain discrete-time linear system of the form
\[
  x_{k+1} = A_k x_k + B_k u_k + C_k w_k,
\]
where $x_k\in\Re^{n_x}$ is a state, $u_k\in\Re^{n_u}$ is a control input, and $w_k\in\Re^{n_w}$ is a disturbance vector.
We assume that the system matrices $A_k\in\Re^{n_x\times n_x}$, $B_k\in\Re^{n_x\times n_u}$, and $C_k\in\Re^{n_x\times n_w}$ are known.

We are interested in controlling the system such that it minimizes the following cost function over a finite prediction horizon $N\in\Nat$:
\begin{equation}\label{eqn:cost_fcn}
  J(x_0,\mbf{u},\mbf{w}) \eqdef \sum_{k=1}^N (x_k^T Q_k x_k + 2 q_k^T x_k) + \sum_{k=0}^{N-1} (u_k^T R_k u_k + 2 r_k^T u_k),
\end{equation}
where $Q_k \succeq 0$ and $R_k \succ 0$ for all $k$, and $\mbf{u}\in\Re^{N_u}$ and $\mbf{w}\in\Re^{N_w}$ denote the entire vectors of control inputs and disturbances, \ie
\begin{align*}
  \mbf{u} &\eqdef (u_0,u_1,\ldots,u_{N-1}) \\
  \mbf{w} &\eqdef (w_0,w_1,\ldots,w_{N-1}),
\end{align*}
with $N_u \eqdef N n_u$ and $N_w \eqdef N n_w$.
Using this notation, the cost function in \eqref{eqn:cost_fcn} can be written in the following compact form:
\[
  J(x_0,\mbf{u},\mbf{w}) = \mbf{w}^T \mbf{C} \mbf{w} + 2 \big( \mbf{c} + \mbf{D}^T \mbf{u} \big)^T \mbf{w} + \mbf{u}^T \mbf{B} \mbf{u} + 2 \mbf{b}^T \mbf{u} + 2 \mbf{a}^T x_0 + x_0^T \mbf{A} x_0 ,
\]
where $\mbf{a}$, $\mbf{A}$, $\mbf{b}$, $\mbf{B}$, $\mbf{c}$, $\mbf{C}$, $\mbf{D}$ are problem data of appropriate dimensions, whose expressions can be found in \cite[\Prop~2]{Bertsimas:2007}; we emphasize that $\mbf{B}\succ 0$ and that only $\mbf{b}$ and $\mbf{c}$ depend on~$x_0$.

Note that \eqref{eqn:cost_fcn} is an uncertain quantity as it depends on the disturbance vector $\mbf{w}$.
Most approaches in stochastic optimization assume $\mbf{w}$ is a random variable and proceed to minimize the expectation of the cost function in \eqref{eqn:cost_fcn}.
An effective way to solve such a problem is via dynamic programming, but this approach is unable to deal tractably with even simple constraints on the control input \cite{Bertsimas:2007}.

\subsection{Robust Control}\label{subsec:lqr:robust}
The authors in \cite{Bertsimas:2007} do not consider any particular distribution for $\mbf{w}$, but assume instead that $\mbf{w}$ belongs to some norm-bounded set
\[
  \mcf{W}_\gamma \eqdef \left\lbrace \mbf{w}\in\Re^{N_w} \mid \norm{\mbf{w}}_2 \le \gamma \right\rbrace,
\]
where $\gamma>0$ can be seen as a bound on the energy of the disturbance signal over the prediction horizon $N$.
We now search for a sequence of control inputs $\mbf{u}^\star$ that solves
\begin{equation}\label{eqn:min-max}
  \min_{\mbf{u}\in\mcf{U}} \, \max_{\mbf{w}\in\mcf{W}_\gamma} \, J(x_0,\mbf{u},\mbf{w}),
\end{equation}
where $\mcf{U}\subseteq\Re^{N_u}$ is a closed convex set.
We can apply this approach in a receding horizon fashion, \ie we compute a solution $\mbf{u}^\star$ of \eqref{eqn:min-max} and implement only the first $n_u$ components to the system, and then reinitialize \eqref{eqn:min-max} to the newly observed state and repeat the procedure.
This means that we need to solve an instance of problem~\eqref{eqn:min-max}, which is parameterized in $x_0$, within the sampling time of the system.

It was shown in \cite[\Thm~3]{Bertsimas:2007} that \eqref{eqn:min-max} can be reformulated as the following SDP:
\begin{equation}\label{eqn:sdp}
  \begin{array}{ll}
    \underset{(\mbf{y},\lambda,z)}{\rm minimize} & z \\
    \text{subject to} & \mbf{y}\in\mcf{Y}, \quad \lambda \ge 0 \\
    & \begin{bmatrix} I & \mbf{y} & \mbf{F} \\ \mbf{y}^T & z-\gamma^2 \lambda & -\mbf{h}^T \\ \mbf{F}^T & -\mbf{h} & \lambda I - \mbf{C} + \mbf{F}^T \mbf{F} \end{bmatrix} \succeq 0,
  \end{array}
\end{equation}
where $\mbf{h}=\mbf{c}-\mbf{D}^T\mbf{B}^{-1}\mbf{b}$, $\mbf{F}=\mbf{B}^{-1/2}\mbf{D}$, $\mbf{u}=\mbf{B}^{-1/2}\mbf{y}-\mbf{B}^{-1}\mbf{b}$, and $\mcf{Y} = \lbrace \mbf{B}^{1/2}\mbf{u} + \mbf{B}^{-1/2}\mbf{b} \mid \mbf{u}\in\mcf{U}\rbrace$.
Although \eqref{eqn:sdp} is a convex optimization problem, it is often not possible to solve it in real time.
This issue motivated the authors in \cite{Bertsimas:2007} to approximate \eqref{eqn:sdp} with an SOCP, which is better suited to real-time applications.

We next show that \eqref{eqn:min-max} can be represented \emph{exactly} as an SOCP, provided that $\mcf{U}$ is SOC-representable.
Using the epigraph reformulation, \eqref{eqn:min-max} is equivalent to
\begin{equation}\label{eqn:epigraph}
  \begin{array}{ll}
    \underset{(\mbf{u},z)}{\rm minimize} & z \\
    \text{subject to} & \mbf{u}\in\mcf{U} \\
    & \mbf{w}^T \mbf{C} \mbf{w} + 2 \big( \mbf{c} + \mbf{D}^T \mbf{u} \big)^T \mbf{w} \le f(x_0,\mbf{u},z), \quad \forall \mbf{w} \colon \mbf{w}^T \mbf{w} \le \gamma^2,
  \end{array}
\end{equation}
where $f(x_0,\mbf{u},z) \eqdef z - \mbf{u}^T \mbf{B} \mbf{u} - 2 \mbf{b}^T \mbf{u} - 2 \mbf{a}^T x_0 - x_0^T \mbf{A} x_0$ is a concave function in $(\mbf{u},z)$.
Let $S\in\Re^{N_w\times N_w}$ be a nonsingular matrix for which
\[
  S^T S = \diag(\sigma_1,\ldots,\sigma_{N_w})
  \qquad\text{and}\qquad
  S^T \mbf{C} S = \diag(\tau_1,\ldots,\tau_{N_w}).
\]
Note that such a matrix exists by Remark~\ref{rem:pos-def}.
Since $\gamma>0$, $\mbf{w}_0=0$ is in the interior of $\mcf{W}_\gamma$.
Thus, we can apply \Lem~\ref{lem:simplified-s}, which states that \eqref{eqn:epigraph} is equivalent to
\begin{equation}\label{eqn:lqr:socp}
  \begin{array}{ll}
    \underset{(\mbf{u},\lambda,t)}{\rm minimize} & \mbf{u}^T \mbf{B} \mbf{u} + 2 \mbf{b}^T \mbf{u} + \boldsymbol{1}^T t + \gamma^2 \lambda \\
    \text{subject to} & \mbf{u}\in\mcf{U}, \quad \lambda \ge 0 \\
    & \Big( \big[ S^T (\mbf{c} + \mbf{D}^T \mbf{u}) \big]_i, t_i, \lambda \sigma_i - \tau_i \Big) \!\in\! \mcf{Q}, \quad \forall i\in[N_w],
  \end{array}
\end{equation}
where we removed the constant terms in the objective function.
If $\mcf{U}$ is polyhedral or SOC-representable, then \eqref{eqn:lqr:socp} can be reformulated as an SOCP, hence belongs to the same complexity class as the problem of minimizing \eqref{eqn:cost_fcn} when $\mbf{w}=0$.
Note that $S$, $\sigma_i$, and $\tau_i$ do not depend on $x_0$ and thus can be precomputed offline.

Observe that the LMI in \eqref{eqn:sdp} has dimension $N_u+N_w+1$, \ie it scales linearly with $N$.
On the other hand, the dimension of each SOC constraint in \eqref{eqn:lqr:socp} is $n_i=2$, but the number of these constraints is equal to $N_w$.

\begin{remark}
  Any feasible solution~$(\mbf{u}',z')$ of \eqref{eqn:epigraph} is robust against the worst-case realization of $\mbf w$ that resides within $\mcf{W}_\gamma$, where a larger value of $\gamma$ promotes the robustness of~$(\mbf{u}',z')$, while the robustness of the solutions often induces probabilistic guarantees.
  For instance, if the disturbances $\mbf{w}$ are perceived as independent and normally distributed random variables, then a probabilistic guarantee for the feasibility of~$(\mbf{u}',z')$ is provided in \cite[\Thm~8]{Bertsimas:2007}.
\end{remark}

\subsection{Robust Regret-Optimal Control}\label{subsec:lqr:regret}
In this section we adopt a different control objective, which is common in the online learning community \cite{Shalev-Shwartz:2012}, and has recently received an increasing interest in the control community as well \cite{Gradu:2020,Goel:2020}.
Instead of minimizing the worst-case cost as in \eqref{eqn:min-max}, we wish to minimize the worst-case \emph{regret}, \ie find a $\mbf{u}^\star$ that solves
\begin{equation}\label{eqn:min-regret}
  \min_{\mbf{u}\in\mcf{U}} \, \max_{\mbf{w}\in\mcf{W}_\gamma} \big( J(x_0,\mbf{u},\mbf{w}) - \min_{\mbf{v}\in\mcf{V}} J(x_0,\mbf{v},\mbf{w}) \big),
\end{equation}
where $\mcf{V}\supseteq\mcf{U}$.
In \eqref{eqn:min-regret} we compare the worst-case performance of the system driven by control input $\mbf{u}$ against using control input $\mbf{v}$, which
(i) knows the whole vector of disturbances $\mbf{w}$ in advance, and
(ii) can choose the sequence of control inputs from a possibly larger set $\mcf{V}$.
Therefore, the optimal value of problem \eqref{eqn:min-regret} is always nonnegative.

We next consider the case where $\mcf{V}=\Re^{N_u}$.
The inner minimization problem in \eqref{eqn:min-regret} reduces to
\[
  \min_{\mbf{v}} \, J(x_0,\mbf{v},\mbf{w}),
\]
which has the  closed-form solution $\mbf{v}^\star=-\mbf{B}\inv(\mbf{b}+\mbf{D}\mbf{w})$.
Plugging $\mbf{v}^\star$ into \eqref{eqn:min-regret}, we obtain
\[
  \min_{\mbf{u}\in\mcf{U}} \, \max_{\mbf{w}\in\mcf{W}_\gamma} \Big( \mbf{w}^T \mbf{D}^T \mbf{B}\inv \mbf{D} \mbf{w} + 2 \big( \mbf{D}^T \mbf{B}\inv \mbf{b} + \mbf{D}^T\mbf{u} \big)^T \mbf{w} + \mbf{u}^T \mbf{B} \mbf{u} + 2\mbf{b}^T\mbf{u} + \mbf{b}^T \mbf{B}\inv \mbf{b} \Big),
\]
which has the same structure as \eqref{eqn:min-max}.
Using similar arguments as in Section~\ref{subsec:lqr:robust}, it follows that \eqref{eqn:min-regret} is equivalent to
\[
  \begin{array}{ll}
    \underset{(\mbf{u},\lambda,t)}{\rm minimize} & \mbf{u}^T \mbf{B} \mbf{u} + 2 \mbf{b}^T \mbf{u} + \boldsymbol{1}^T t + \gamma^2 \lambda \\
    \text{subject to} & \mbf{u}\in\mcf{U}, \quad \lambda \ge 0 \\
    & \Big(  \big[ S^T \mbf{D}^T  (\mbf{B}\inv \mbf{b} + \mbf{u}) \big]_i, t_i, \lambda \sigma_i - \tau_i  \Big) \in \mcf{Q}, \quad \forall i\in[N_w],
  \end{array}
\]
where $S\in\Re^{N_w\times N_w}$ is a nonsingular matrix for which
\[
  S^T S = \diag(\sigma_1,\ldots,\sigma_{N_w})
  \qquad\text{and}\qquad
  S^T \mbf{D}^T\mbf{B}\inv\mbf{D} S = \diag(\tau_1,\ldots,\tau_{N_w}).
\]

\subsection{Distributionally Robust Control}
Robust control may be too conservative as no statistical knowledge about the disturbance $\mbf{w}$ is incorporated in \eqref{eqn:min-max}.
One popular remedy is to adopt a distributionally robust approach \cite{Wiesemann:2014}, where one seeks an optimal control input vector $\mbf{u}$ that minimizes the expected value of the cost function $J$ with respect to the worst-case probability distribution $\mbb{P}$ that resides within the ambiguity set $\mathscr{P}$, \ie
\begin{equation}\label{eqn:min-sup-exp}
  \min_{\mbf{u}\in\mcf{U}} \, \sup_{\mbb{P}\in\mathscr{P}} \, \mbb{E}_\mbb{P} [ J(x_0,\mbf{u},\tilde{\mbf{w}})],
\end{equation}
where $\mathscr{P} = \{ \mbb{P} \in \mathscr{P}_0 (\mcf{W}_\gamma) \mid \mbb{E}_\mbb{P} [H \tilde{\mbf{w}}] \le \mu \}$ contains all probability distributions that satisfy $m$ first-order moment conditions that are modeled through $H\in\Re^{m \times N_w}$ and $\mu\in\Re^{m}$, and $\mathscr{P}_0 (\mcf{W}_\gamma)$ represents a family of all distributions supported on $\mcf{W}_\gamma$.
The considered ambiguity set $\mathscr{P}$ can incorporate, for instance, the mean and mean-absolute deviation of the disturbances; we refer the interested reader to \cite{Wiesemann:2014,Chen:2020} for more details on the expressiveness of $\mathscr{P}$.

Thanks to \cite[\Prop~3.4]{Shapiro:2001}, we can reformulate~\eqref{eqn:min-sup-exp} as the following semi-infinite problem under mild regularity conditions:
\begin{equation}\label{eqn:min-max-exp}
  \begin{array}{llr}
    \underset{(\mbf{u},\alpha,\beta)}{\rm minimize} & \alpha + \mu^T \beta \\
    \text{subject to} & \mbf{u}\in\mcf{U}, \quad \beta \ge 0 \\
    & \alpha + \mbf{w}^T H^T \beta \ge  J(x_0,\mbf{u},\mbf{w}), \quad \forall \mbf{w}\in\mcf{W}_\gamma.
  \end{array}
\end{equation}
Since \eqref{eqn:min-max-exp} has the same structure as~\eqref{eqn:epigraph}, we can follow similar arguments as in Section~\ref{subsec:lqr:robust} to show that it is equivalent to
\[
  \begin{array}{ll}
    \underset{(\mbf{u},\beta,\lambda,t)}{\rm minimize} & \mbf{u}^T \mbf{B} \mbf{u} + 2 \mbf{b}^T \mbf{u} + \boldsymbol{1}^T t + \gamma^2 \lambda + \mu^T \beta \\
    \text{subject to} & \mbf{u}\in\mcf{U}, \quad \beta \ge 0, \quad \lambda \ge 0 \\
    & \Big( \big[ S^T (\mbf{c} + \mbf{D}^T \mbf{u} - \tfrac{1}{2} H^T \beta ) \big]_i, t_i, \lambda \sigma_i - \tau_i \Big) \in \mcf{Q}, \quad \forall i\in[N_w],
  \end{array}
\]
where $S$, $\sigma_i$, and $\tau_i$ are the same as in Section~\ref{subsec:lqr:robust}.

Similarly, the following distributionally robust version of the regret-optimal control problem \eqref{eqn:min-regret} with $\mcf{V}= \Re^{N_u}$
\begin{equation*}
  \min_{\mbf{u}\in\mcf{U}} \, \sup_{\mbb{P}\in\mathscr{P}} \mbb{E}_\mbb{P} \big[ J(x_0,\mbf{u},\mbf{w}) - \min_{\mbf{v}} J(x_0,\mbf{v},\mbf{w}) \big],
\end{equation*}
can be analogously reformulated as
\[
  \begin{array}{ll}
    \underset{(\mbf{u},\beta,\lambda,t)}{\rm minimize} & \mbf{u}^T \mbf{B} \mbf{u} + 2 \mbf{b}^T \mbf{u} + \boldsymbol{1}^T t + \gamma^2 \lambda +\mu^T \beta \\
    \text{subject to} & \mbf{u}\in\mcf{U}, \quad \beta \ge 0,  \quad \lambda \ge 0 \\
    & \Big( \big[ S^T \big( \mbf{D}^T  (\mbf{B}\inv \mbf{b} + \mbf{u}) - \frac{1}{2} H^T\beta \big) \big]_i, t_i, \lambda \sigma_i - \tau_i  \Big) \in \mcf{Q}, \quad \forall i\in[N_w].
  \end{array}
\]

\section{Reconfigurable Terminal Constraints in MPC}\label{sec:mpc}
For a linear time-invariant discrete-time system
\begin{equation}\label{eqn:sys}
  x_{k+1} = A x_k + B u_k,
\end{equation}
with $A\in\Re^{n_x\times n_x}$, $B\in\Re^{n_x\times n_u}$, state constraints $x_k\in\mcf{X}$, and input constraints $u_k\in\mcf{U}$, the MPC law is defined through the solution of a constrained finite-time optimal control problem in a receding horizon fashion.
A well-known approach for ensuring asymptotic stability and recursive feasibility of the MPC scheme is to include a particular terminal cost and/or terminal constraints in the problem \cite{Sznaier:1987,Mayne:2000}.

The MPC optimization problem with these terminal ingredients can be formulated as
\begin{equation}\label{eqn:mpc}
  \begin{array}{ll}
    \underset{(\mbf{u},\mbf{x})}{\rm minimize} & \displaystyle \sum_{k=0}^{N-1} l(x_k,u_k) + \psi(x_N) \\
    \text{subject to} & x_0 = x_{\rm init} \\
    & x_{k+1} = f(x_k,u_k), \hspace{0.2cm} \forall k\in\{0,\ldots,N-1\} \\
    & x_k \in \mcf{X}, \hspace{1.85cm} \forall k\in\{1,\ldots,N-1\} \\
    & u_k \in \mcf{U}, \hspace{1.9cm} \forall k\in\{0,\ldots,N-1\} \\
    & x_N \in \mcf{T},
  \end{array}
\end{equation}
where $l\colon\mcf{X}\times\mcf{U}\to\Re$ is the stage cost, $\mcf{T}\subseteq\mcf{X}$ the terminal set, $\psi\colon\mcf{T}\to\Re$ the terminal cost, $x_{\rm init}\in\mcf{X}$ the initial state, and
\begin{align*}
  \mbf{u} &\eqdef (u_0,u_1,\ldots,u_{N-1}) \\
  \mbf{x} &\eqdef (x_0,x_1,\ldots,x_N).
\end{align*}
While $l$, $\mcf{X}$, and $\mcf{U}$ are typically part of control specifications, the terminal ingredients $\psi$ and $\mcf{T}$ are design choices and can be used to ensure asymptotic stability and recursive feasibility of the MPC scheme.
Although the terminal ingredients are usually computed offline, as part of the controller design phase, the authors in \cite{Simon:2014} consider a \emph{polyhedral} terminal set, which can be scaled and translated online, \ie its center and the scaling factor are computed within the MPC optimization problem.
The flexibility in choosing the terminal set online results in an enlarged feasible region of problem \eqref{eqn:mpc}, which in general improves the performance of the controller.

We consider polyhedral state and input constraints, and an \emph{ellipsoidal} terminal set given by
\begin{subequations}\label{eqn:sets}
\begin{align}
  \mcf{X} &= \lbrace x\in\Re^{n_x} \mid Ex \le f \rbrace \\
  \mcf{U} &= \lbrace u\in\Re^{n_u} \mid Gu \le h \rbrace \\
  \mcf{T} &= \lbrace x\in\Re^{n_x} \mid (x-c)^T P (x-c) \le r^2 \rbrace, \label{eqn:sets:T}
\end{align}
\end{subequations}
where $E\in\Re^{n\times n_x}$, $f\in\Re^n$, $G\in\Re^{m\times n_u}$, $h\in\Re^m$, and $P\in\symm_{++}^{n_x}$ are assumed to be known.
While most approaches assume $c\in\Re^{n_x}$ and $r>0$, that determine the center and scaling of the terminal set, are fixed quantities, we consider them as decision variables.
A similar form of the reconfigurable terminal set is used in \cite{Aboudonia:2020} in the context of distributed MPC.

\subsection{Online Computation of the Terminal Set}\label{subsec:mpc:online}
A standard way to ensure asymptotic stability and recursive feasibility of the MPC scheme is to assume that after applying the sequence of $N$ control inputs computed in \eqref{eqn:mpc}, a linear terminal controller $\kappa\colon\mcf{T}\to\mcf{U}\colon x\mapsto Kx$, with $K\in\Re^{n_u\times n_x}$, takes over controlling the system.
We assume that the terminal controller is stabilizing so that the terminal closed-loop dynamics $x_{k+1}=\Acl x_k$, with $\Acl \eqdef A+BK$, is stable.

The closed-loop system controlled with the MPC scheme is asymptotically stable and recursively feasible if
(i) $\mcf{T}$ is a positively invariant set of system \eqref{eqn:sys} controlled with a terminal controller $\kappa$,
(ii) $\mcf{T}$ is included in $\mcf{X}$, and
(iii) running the terminal controller on $\mcf{T}$ does not violate the input constraints \cite{Mayne:2000}.
In other words, we want to ensure that for all $x\in\mcf{T}$, the following inclusions hold:
\begin{subequations}\label{eqn:pos-inv}
\begin{align}
  \Acl x &\in \mcf{T} \label{eqn:pos-inv:T} \\
  x &\in \mcf{X} \label{eqn:pos-inv:X} \\
  Kx &\in \mcf{U}. \label{eqn:pos-inv:U}
\end{align}
\end{subequations}
The authors in \cite{Aboudonia:2020} and \cite{Darivianakis:2020} consider similar requirements and reformulate the robust constraints in \eqref{eqn:pos-inv} as LMIs.
We next show how to reformulate them as linear and SOC constraints.

\subsubsection{Positive invariance}
The robust constraint in \eqref{eqn:pos-inv:T}, with $\mcf{T}$ given by \eqref{eqn:sets:T}, can be written as
\[
  (\Acl x-c)^T P (\Acl x-c) \le r^2, \quad \forall x \colon (x-c)^T P (x-c) \le r^2.
\]
We rewrite the constraint above in terms of auxiliary variables $s=r\inv (x-c)$ and $\hat{c}=r\inv c$,
\begin{equation}\label{eqn:pos-inv-semiinf}
  s^T \Acl^T P \Acl s + 2 \big( \Acl^T P (\Acl-I)\hat{c} \big)^T s \le 1 - \hat{c}^T M \hat{c}, \quad \forall s \colon s^T P s \le 1,
\end{equation}
where $M \coloneqq (\Acl\!-\!I)^T P (\Acl\!-\!I) \succeq 0$.
Let $S\in\Re^{n_x\times n_x}$ be a nonsingular matrix for which
\begin{equation}\label{eqn:mpc-S}
  S^T P S = \diag(\pi_1,\ldots,\pi_{n_x})
  \qquad\text{and}\qquad
  S^T \Acl^T P \Acl S = \diag(\alpha_1,\ldots,\alpha_{n_x}).
\end{equation}
Since $s_0^T P s_0 < 1$ for $s_0=0$, we can apply \Lem~\ref{lem:simplified-s}, which states that \eqref{eqn:pos-inv-semiinf} holds if and only if there exists $(\lambda,t)\in\Re_+\times\Re^{n_x}$ such that
\begin{align*}
  & \boldsymbol{1}^T t + \lambda \le 1-\hat{c}^T M \hat{c} \\
  & \Big( \big[ S^T \Acl^T P (\Acl-I)\hat{c} \big]_i, t_i, \lambda \pi_i - \alpha_i \Big) \in \mcf{Q}, \quad \forall i\in[n_x].
\end{align*}
Multiplying both constraints by $r$ and using \Lem~\ref{lem:hyperbolic}, we obtain
\begin{equation}\label{eqn:pos-inv-final}
  \begin{split}
    & \big\| \big( 2M^{1/2}c, \boldsymbol{1}^T \hat{t} + \hat{\lambda} \big) \big\|_2 \le 2r - \boldsymbol{1}^T \hat{t} - \hat{\lambda} \\
    & \Big( \big[ S^T \Acl^T P (\Acl-I) c \big]_i, \hat{t}_i, \hat{\lambda} \pi_i - r \alpha_i \Big) \in \mcf{Q}, \quad \forall i\in[n_x],
  \end{split}
\end{equation}
where $\hat{\lambda}=r\lambda$ and $\hat{t}=rt$.

\subsubsection{State and input constraints}
The robust constraints in \eqref{eqn:pos-inv:X} and \eqref{eqn:pos-inv:U}, with $\mcf{X}$, $\mcf{U}$, and $\mcf{T}$, given by \eqref{eqn:sets}, can be written as
\begin{subequations}\label{eqn:robust}
\begin{align}
  E x \le f, \quad & \forall x \colon (x-c)^T P (x-c) \le r^2 \label{eqn:robust:state} \\
  GK x \le h, \quad & \forall x \colon (x-c)^T P (x-c) \le r^2. \label{eqn:robust:input}
\end{align}
\end{subequations}
Even if $\mcf{X}$ and $\mcf{U}$ contained quadratic constraints, we would be able to reformulate \eqref{eqn:pos-inv:X} and \eqref{eqn:pos-inv:U} as a set of linear and SOC constraints using similar arguments as in the previous section.
As $\mcf{X}$ and $\mcf{U}$ are polyhedral, we next show that \eqref{eqn:robust} can be reformulated using only linear constraints.

Introducing an auxiliary variable $p=P^{1/2}(x-c)$, we can rewrite \eqref{eqn:robust:state} as
\[
  Ec + E P^{-1/2} p \le f, \quad \forall p \colon \norm{p}_2 \le r,
\]
which can be reformulated in a similar fashion as done in \cite[\S 4.3.1]{Boyd:2004}.
In particular, we consider the $j$-th constraint
\[
  e_j^T c + \hat{e}_j^T p \le f_j, \quad \forall p \colon \norm{p}_2 \le r,
\]
where $e_j^T$ and $\hat{e}_j^T$ denote the $j$-th rows of $E$ and $E P^{-1/2}$, respectively.
This robust constraint is satisfied if and only if
\[
  e_j^T c + \sup_{\norm{p}_2\le r} \! \big( \hat{e}_j^T p \big) \le f_j,
\]
which is equivalent to 
\[
  e_j^T c + \norm{\hat{e}_j}_2 r \le f_j.
\]
Therefore, \eqref{eqn:robust:state} holds if and only if the inequality above holds for all $j\in[n]$, \ie
\begin{equation}\label{eqn:state-constr-final}
  e_j^T c + \norm{\hat{e}_j}_2 r \le f_j, \quad j\in[n].
\end{equation}
Similarly, \eqref{eqn:robust:input} is equivalent to
\begin{equation}\label{eqn:input-constr-final}
  g_j^T c + \norm{\hat{g}_j}_2 r \le h_j, \quad j\in[m],
\end{equation}
where $g_j^T$ and $\hat{g}_j^T$ denote the $j$-th rows of $GK$ and $GK P^{-1/2}$, respectively.

\subsection{Tractable Reformulation of the MPC Problem}
The constraint $x_N\in\mcf{T}$ is equivalent to
\[
  \| P^{1/2} (x_N-c) \|_2 \le r.
\]
Therefore, problem~\eqref{eqn:mpc} with $\mcf{X}$, $\mcf{U}$, and $\mcf{T}$ given by \eqref{eqn:sets}, together with the positive invariance conditions in \eqref{eqn:pos-inv}, is equivalent to
\[
  \begin{array}{ll}
    \underset{(\mbf{u},\mbf{x},c,r,\hat{\lambda}, \hat{t})}{\rm minimize} & \displaystyle \sum_{k=0}^{N-1} l(x_k,u_k) + \psi(x_N) \\
    \text{subject to} & x_0 = x_{\rm init} \\
    & x_{k+1} = A x_k + B u_k, \hspace{0.4cm} \forall k\in\{0,\ldots,N-1\} \\
    & E x_k \le f, \hspace{2.2cm} \forall k\in\{1,\ldots,N-1\} \\
    & G u_k \le h, \hspace{2.2cm} \forall k\in\{0,\ldots,N-1\} \\
    & \| P^{1/2} (x_N-c) \|_2 \le r \\
    & \eqref{eqn:pos-inv-final}, \, \eqref{eqn:state-constr-final}, \, \eqref{eqn:input-constr-final}.
  \end{array}
\]
In the standard case, in which $l$ and $\psi$ are convex quadratic functions as in \eqref{eqn:cost_fcn}, the problem above can be represented as an SOCP.

\section{Numerical Example}\label{sec:numerics}
To demonstrate the computational benefit of using a simplified \SL-lemma, we evaluate the performance of the respective SDP and SOCP reformulations~\eqref{eqn:sdp} and~\eqref{eqn:lqr:socp}
of robust LQC problem \eqref{eqn:min-max} given in Section~\ref{subsec:lqr:robust}.
We adapt a numerical example given in \cite{Bertsimas:2007} and consider a simple linear system with time-invariant state, control, and disturbance matrices $A_k=B_k=C_k=1$ for all $k$ and initial state $x_0=-1$.
The cost function is given by $Q_k = R_k = 0.9^k$ and $q_k=r_k=0$ for all $k\in\{0,1,\ldots,N\}$.
We set $\gamma=0.1$ and define the input constraint set as $\mcf{U} = \lbrace \mbf{u}\in\Re^{N_u} \mid -0.4 \cdot \boldsymbol{1} \le \mbf{u} \le 0.4 \cdot \boldsymbol{1} \rbrace$.

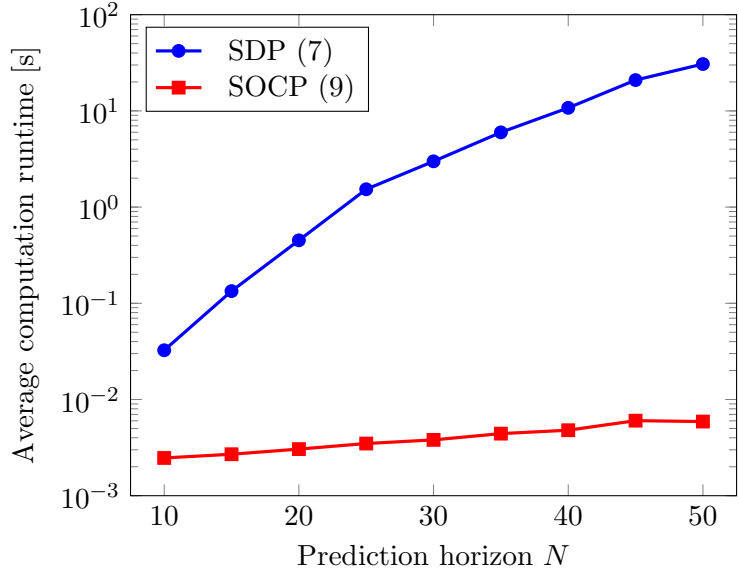
\begin{figure}[t]
  \centering
  \small
  \begin{tikzpicture}
    \begin{axis}[
      width=\mywidth,
      height=\myheight,
      ymin=1e-3,
      ymax=1e2,
      xmin=7.5,
      xmax=52.5,
      xlabel={Prediction horizon $N$},
      ylabel={Average computation runtime},
      y unit=\si{\second},
      yticklabel pos=left,
      ymode=log,
      legend cell align=left,
      legend style={at={(.025,.975)}, nodes={scale=1}, anchor=north west, fill=white, fill opacity=1, draw opacity=1,text opacity=1},
      every axis plot/.append style={very thick}
      ]
      \addlegendimage{blue, mark=*}
      \addlegendentry{\small SDP \eqref{eqn:sdp}};
      \addlegendimage{red, mark=square*}
      \addlegendentry{\small SOCP \eqref{eqn:lqr:socp}};
      \addplot [blue, mark=*, mark options={scale=1, fill=blue}]
      table[x=N, y=SDP, col sep=comma] {runtime_mean.csv};
      \addplot [red, mark=square*, mark options={scale=1, fill=red}]
      table[x=N, y=SOCP, col sep=comma] {runtime_mean.csv};
    \end{axis}
  \end{tikzpicture}
  \caption{Average computation runtimes of SDP and SOCP reformulations of the robust LQC problem given in \eqref{eqn:min-max}.}
  \label{fig:numerics}
\end{figure}

We perform numerical tests on a Linux-based system with an i9-9900K @ 3.6 GHz (8 cores) processor and 64 GB of DDR4 3200Mhz RAM, and use CVXPY \cite{Diamond:2016} and MOSEK \cite{mosek} to solve the optimization problems.
Figure~\ref{fig:numerics} shows the average computation times required to solve the resulting problems over $10$ runs for each value of the prediction horizon $N$.
It can be seen that for longer prediction horizons solving the SOCP is more than 1000 times faster.

\section{Conclusion}\label{sec:conclusion}
We have revisited the so-called simplified \SL-lemma and have shown that various optimization problems arising in control, that are usually formulated as SDPs, can be represented as SOCPs.
This is particularly important for real-time applications that require the repeated solution of such an optimization problem, as SOCPs can be solved faster and require simpler numerical solvers.
Our results pave the way for similar developments in other control schemes that currently use the standard \SL-lemma.

\section*{Acknowledgements}
We are grateful to Ahmed Aboudonia for helpful discussions on reconfigurable terminal constraints in MPC.
\ackERC

\bibliography{refs}

\end{document}

%% file: defs.tex
\usepackage{parskip}
\usepackage{fullpage}
\usepackage{abstract}

\usepackage{titlesec}
\def\headingfont{\sffamily\bfseries\boldmath}
\titleformat*{\section}{\headingfont\Large}
\titleformat*{\subsection}{\headingfont\large}
\titleformat*{\subsubsection}{\headingfont\normalsize}
\titleformat*{\paragraph}{\headingfont\normalsize}
\titleformat*{\subparagraph}{\headingfont\normalsize}

\usepackage{amsmath,amssymb,mathtools,amsfonts,amsthm,mathrsfs}
\usepackage[english]{babel}
\usepackage[titletoc,title]{appendix}
\usepackage{enumitem}
\setlist{nolistsep}
\usepackage{hyperref}
\hypersetup{%
  colorlinks=true,
  linktoc=all,
  linkcolor=black,
  citecolor= black,
  urlcolor=blue
}

\newtheoremstyle{exampstyle}
  {.5\baselineskip}
  {\topsep}
  {}
  {}
  {\bfseries}
  {.}
  {.5em}
  {}
\theoremstyle{exampstyle}

\newtheorem{theorem}{Theorem}[section]

\newtheorem{lemma}[theorem]{Lemma}

\newtheorem{remark}[theorem]{Remark}

\newcommand{\ie}{\emph{i.e.},\ }

\newcommand{\mcf}{\mathcal}
\newcommand{\mbf}[1]{\boldsymbol{#1}}
\newcommand{\mbb}{\mathbb}

\renewcommand{\Re}{\mbb{R}}
\newcommand{\Nat}{\mbb{N}}
\newcommand{\symm}{\mbb{S}}

\newcommand{\Prop}{Prop.}

\newcommand{\Thm}{Thm.}
\newcommand{\Lem}{Lemma}

\DeclareMathOperator{\diag}{diag}                                                 

\newcommand{\eqdef}{\coloneqq}
\newcommand{\inv}{^{-1}}
\newcommand{\norm}[1]{\lVert#1\rVert}
\newcommand{\SL}{$\mcf{S}$}
\newcommand{\Acl}{A_{\rm cl}}

\usepackage[binary-units=true]{siunitx}
\usepackage{tikz}
\usepackage{pgfplots}
\pgfplotsset{compat=newest}
\usepgfplotslibrary{units}
\pgfplotsset{
    /pgfplots/ybar legend/.style={
    /pgfplots/legend image code/.code={%
       \draw[##1,/tikz/.cd,
       yshift=-0.3em
       ]
        (0cm,0cm) rectangle (5pt,0.8em);},
   }
}

\newlength{\mywidth}
\newlength{\myheight}